\documentclass[twocolumn,letterpaper,10pt]{article}
\usepackage{titlesec}
\usepackage{pslatex}
\usepackage{url}
\usepackage{indentfirst}
\usepackage{balance}
\usepackage{parskip}
\usepackage{enumitem}
\AtBeginDocument{\RequirePackage{lmodern, times}}
\usepackage{amsmath,amsfonts,amsthm}
\usepackage{algorithmic}
\usepackage{algorithm}
\usepackage[table]{xcolor}
\usepackage[pdftex]{graphicx}

\titleformat{\section}{\normalfont\fontsize{12}{15}\bfseries}{\thesection.}{0.5em}{}
\titleformat{\subsection}{\normalfont\bfseries}{\thesubsection.}{3pt}{}
\titleformat{\subsubsection}[runin]{\normalfont\bfseries}{\thesubsubsection.}{3pt}{}

\titlespacing{\section}{0pt}{3ex}{3ex}
\titlespacing{\subsection}{0pt}{2ex}{2ex}

\newtheorem{theorem}{Theorem}
\newtheorem{definition}{Definition}
\newtheorem{lemma}{Lemma}

\newtheorem{corollary}{Corollary}

\newtheorem{assumption}{Assumption}

\newcommand{\Nset}{\mathcal{N}}

\newcommand{\Eset}{\mathcal{E}}
\newcommand{\tp}{\textsc{3-Partition}}

\newcommand{\abs}[1]{\lvert #1 \rvert}

\newcommand{\etal}{{\em et~al.~}}

\definecolor{pennblue}{cmyk}{1,0.65,0,0.30}
\definecolor{pennred}{cmyk}{0,1,0.65,0.34}
\definecolor{mygreen}{rgb}{0.10,0.50,0.10}

\begin{document}
\title{
{\fontsize{14}{10}\bfseries A Submodular Approach for\\Electricity Distribution Network Reconfiguration}}



\author{
   Ali Khodabakhsh$^{*}$, Ger Yang$^{*}$, Soumya Basu$^{*}$, Evdokia Nikolova$^{*}$, Michael C. Caramanis$^{\dag}$,\\ Thanasis Lianeas$^{*}$, Emmanouil Pountourakis$^{*}$ \\
   \normalsize $^{*}$ {Department of Electrical and Computer Engineering, University of Texas at Austin}  \\[-3pt]
    \normalsize \{ali.kh, geryang, basusoumya, thanasis, manolis\}@utexas.edu, nikolova@austin.utexas.edu\\
    \normalsize $^{\dag}$ {Department of Mechanical Engineering, Boston University} \\[-3pt]
    \normalsize mcaraman@bu.edu \\[-3pt]
}

\date{} 

\maketitle
\thispagestyle{empty}

\centerline{{\bf\large Abstract}}

{\it Distribution network reconfiguration (DNR) is a tool used by operators to balance line load flows and mitigate losses. As distributed generation and flexible load adoption increases, the impact of DNR on the security, efficiency, and reliability of the grid will increase as well.  Today, heuristic-based actions like branch exchange are routinely taken, with no theoretical guarantee of their optimality. 
This paper considers loss minimization via DNR, which changes the on/off status of switches in the network.  The goal is to ensure a radial final configuration (called a spanning tree in the algorithms literature) that spans all network buses and connects them to the substation (called the root of the tree) through a single path.  We prove that the associated combinatorial optimization problem is strongly NP-hard and thus likely cannot be solved efficiently. We formulate the loss minimization problem as a supermodular function minimization under a single matroid basis constraint, and use existing algorithms to propose a polynomial time local search algorithm for the DNR problem at hand and derive performance bounds.  We show that our algorithm is equivalent to the extensively used branch exchange algorithm, for which, to the best of our knowledge, we pioneer in proposing a theoretical performance bound. Finally, we use a 33-bus network to compare our algorithm's performance to several algorithms published in the literature.}

\section{Introduction}
Distribution networks are usually built as interconnected mesh networks, but are normally configured (via switches) and operated as radial networks (i.e. trees, in graph theoretic terms), to simplify overload protection \cite{pso}. 
The entire network can be thought of as a forest consisting of rooted trees. 
Each tree consists of a substation (root) and a number of customers (users) that are serviced via so-called distribution feeders (distribution lines starting at the substation).  
Switches located throughout the network allow dynamic reconfiguration of the distribution network through switching operations; the opening or closing of a switch corresponds to the removal or addition of an edge, respectively. 

The goal of distribution networks is to deliver the power from substations to users, but notably, substantial losses of up to 13\% occur as electric power flows over distribution lines \cite{sarfi1994survey}. As a result, \emph{Distribution Network Reconfiguration} (DNR) is a major tool focusing on the dynamic identification of a spanning tree  that optimizes a performance measure such as load flow balancing or total line loss minimization. We select the latter, namely the minimization of losses for a given hourly load flow, as the objective of the reconfiguration problem. Similar issues in meshed transmission networks have been addressed in the literature recently (see \cite{goldis2015ac} and references therein).

\noindent{\bf Our results:} In this paper, we analyze the DNR problem via a submodular optimization approach. In particular, we give the following results:
\begin{enumerate}[leftmargin=*]
\item We prove that the DNR problem is strongly NP-hard. We do this through a polynomial reduction from \tp\ problem, which is defined in Section \ref{sec:hardness} (see \cite{gary1979computers} for more details). To the best of our knowledge, the computational hardness of this problem has not been studied so far.
\item We formulate the DNR problem as a supermodular minimization problem subject to a single matroid basis constraint (we define supermodularity and matroid later in Section \ref{subsec:submodularity}). Supermodularity is motivated by the fact that losses are quadratic in the current flowing over each branch of the distribution network. Furthermore, the matroid basis constraint ensures the radial structure and guarantees that all the buses are connected to the substation.

\item We observe that the local search algorithm for solving the supermodular minimization problem is equivalent to the well-known \emph{branch exchange} algorithm. Hence, we obtain the first theoretical result on why the branch exchange algorithm performs well in practice.
\end{enumerate}

The proposed submodular framework sheds some light on the algorithmic structure of the optimization problems in distribution networks. Although for the DNR problem we are mostly providing a theoretical justification for an existing heuristic, as it is evident in other lines of work in energy systems (see \cite{ex1,ex2,ex3} for example), the theoretical study of such problems can help to either find new algorithms or improve the efficiency of existing ones.   

The rest of this paper is organized as follows. Section~\ref{sec:related} reviews related work. Section~\ref{sec:model} gives a concise formulation of the problem; and its computational complexity is studied in Section~\ref{sec:hardness}. The submodular framework is proposed in Section~\ref{sec:supermodular}. Section~\ref{sec:alg} describes the algorithm and its performance guarantee. Section~\ref{sec:simulation}, provides numerical results and comparison with different algorithms. Finally, Section~\ref{sec:conclusion} concludes the work.
\section{Related Work} \label{sec:related}
\begin{figure}[t] 
   \centering
   \includegraphics[width=2.85in]{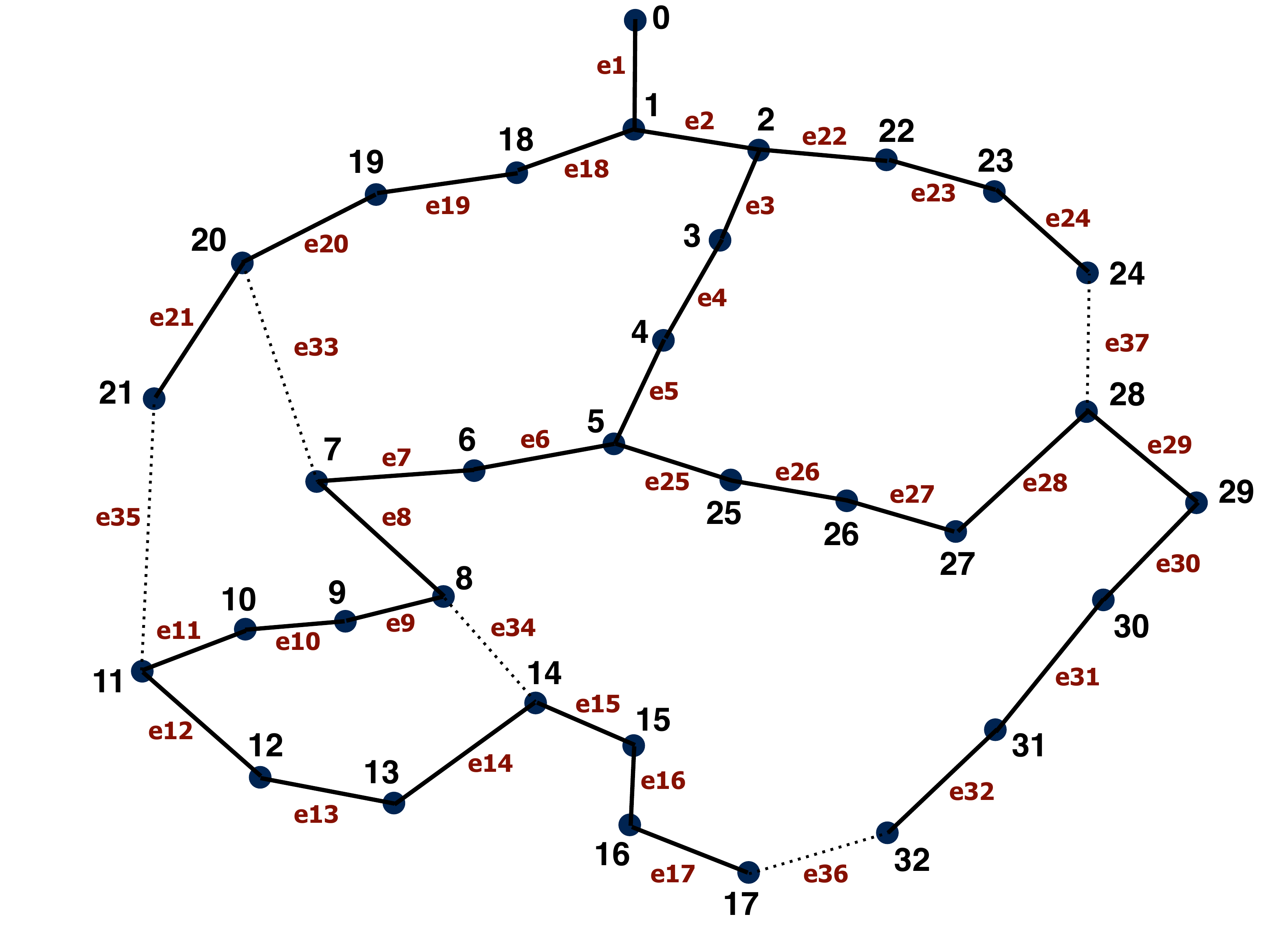} 
   \vspace{-0.1in}
   \caption{33-bus network \cite{khodr2009integral}.}
   \label{fig:33bus}
\end{figure}

DNR has been studied extensively in the literature. One of the most common heuristic algorithms is the \emph{branch exchange} suggested by Civanlar~\etal \cite{civanlar1988distribution} and implemented by Baran and Wu \cite{baran1989network}, who considered loss minimization and load balancing objectives. Starting from a feasible tree configuration, the branch exchange algorithm transfers some loads in each iteration by (i) closing an open switch to create a loop in the network, followed by (ii) opening one of the closed switches in that loop to arrive at another feasible solution with a lower cost. The algorithm terminates when no further improvements are possible. This algorithm has been used as a benchmark against different DNR algorithms with the $12.6$kV network of Fig.~\ref{fig:33bus} employed for numerical comparisons.

An improved branch exchange algorithm was proposed by Miguez \etal \cite{miguez2002improved} who tried to expand the space of available changes in the local search, hence eliminating some local minima of the standard algorithm. The idea of improved branch exchange is to investigate improvement from a pair of exchanges, once there is no improvement by a single branch exchange. Peng and Low \cite{peng2013optimal} proposed an algorithm to do each step of branch exchange efficiently by solving only 3 optimal power flow equations (OPF), regardless of the size of the network. Their algorithm helps to find the best switch to open in order to minimize any convex increasing cost function, assuming that an open switch has already been closed. These improvements still provide no theoretical guarantee on the output of the branch exchange algorithm.

Unlike the branch exchange algorithm that maintains a tree structure during its execution, there are other heuristic algorithms that start with the meshed network (obtained by closing all the tie switches) or the disconnected network (obtained by opening all the switches) and proceed to open/close switches one by one until a radial configuration is achieved \cite{gomes2005new,shirmo1989reconfig,mcdermott1999heuristic}. Shirmohammadi and Hong \cite{shirmo1989reconfig} proposed one such algorithm that starts with the meshed network and proceeds with iterations that open the switch with the smallest current. No theoretical performance guarantees have been obtained for this algorithm.

For small networks such as the 33-bus example of Fig.~\ref{fig:33bus}, the global optimal configuration can be discovered by brute-force enumeration. An efficient enumeration approach proposed in \cite{morton2000efficient}, lists all the spanning trees in a clever way that generates each tree exactly once, and calculates losses by adjusting the losses of the previous spanning tree. The drawback of this method is that it is not practical for larger networks, since a network has exponentially many spanning trees \cite{kocay2016graphs}.

The joint DNR and OPF problem was considered in \cite{khodr2009integral} using Benders decomposition to decompose the global problem to master and slave subproblems. The master level determines the binary variables by solving a mixed-integer non-linear program using CPLEX. The slave level solves the OPF non-linear program using the CONOPT solver. Again, solving integer programs is computationally intractable for large networks.

Many other approaches like genetic algorithms \cite{lin2000distribution,Enacheanu2008radial}, particle swarm optimization \cite{pso}, ant colony algorithms \cite{wu2010study}, artificial neural networks \cite{kim1993artificial,kashem1998artificial}, etc. have been utilized to solve this problem.
A survey of different algorithms for the DNR problem can be found in \cite{sarfi1994survey}. What is conspicuously missing in all these previous works is a rigorous theoretical performance guarantee.

To close this gap, we consider a submodular approach to the DNR problem. Since switching binary decisions render DNR a non-linear combinatorial optimization problem, additional structure like submodularity or supermodularity enables finding an approximate solution efficiently.
\section{Problem Formulation} \label{sec:model}
In this section we present the power flow equations and employ some simplifying assumptions to model the problem in graph theoretic terms. 
We model the distribution network as a graph $\mathcal{G}(\Nset,\Eset)$, where $\Nset$ is the set of buses (nodes) and $\Eset$ is the set of lines (undirected edges). We assume that a single  substation is located at node 0, and the other nodes are load buses with given active and reactive power demands $(p_i,q_i)$, for all $i\in \Nset\backslash \{0\}$. We are looking for a spanning tree rooted at bus $0$ (i.e., a tree that connects all the loads to the root through a single path) which minimizes the total resistive loss.

\begin{figure}[t] 
   \centering
   \includegraphics[width=1.6in]{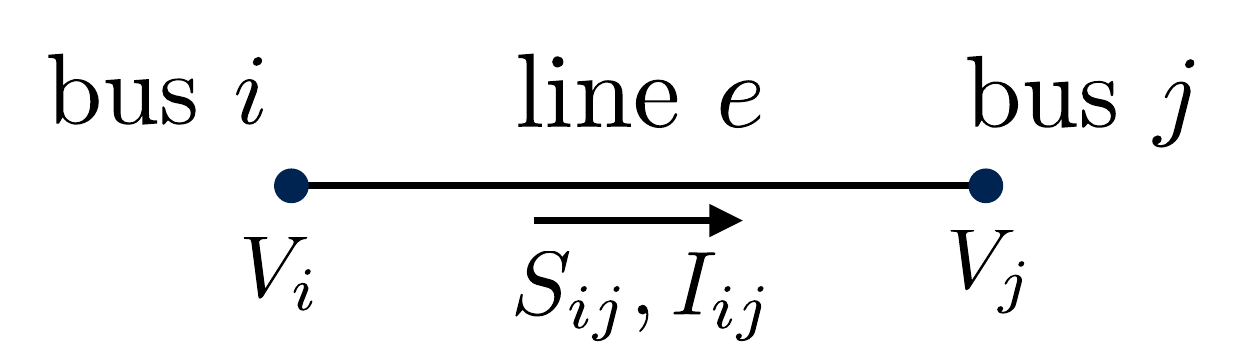} 
   \caption{Power flow variables.}
   \label{fig:flow}
\end{figure}

Letting $V_i=\abs{V_i}e^{\mathbf{i}\theta_i}$ represent the complex voltage at bus $i$, we adopt the relaxed branch model of \cite{peng2013optimal,wu1989relaxed} that allows us to ignore the phase angles of voltages and currents in radial networks. Let $Z_e=R_e+\mathbf{i}X_e$ be the impedance of line $e\in \Eset$. We also use $S_{ij}=P_{ij}+\mathbf{i}Q_{ij}$ to express the branch power flow from bus $i$ to bus $j$, and $I_{ij}$ to express the current from bus $i$ to $j$. A summary of our notation is depicted in Fig.~\ref{fig:flow}.
\vspace{0.1cm}
\begin{assumption}[{\cite[A2]{peng2013optimal}}]
\label{assumption:drop}
Voltage variation across the distribution network can be neglected. Using per unit (p.u.) representation, we assume that $\abs{V_i}=1$ p.u. for all nodes $i\in \Nset$.
\end{assumption}
This assumption is realistic since in practice voltage at every bus is kept within an allowable range such as $(0.95,1.05)$ p.u., and impacts losses (the objective function of DNR) at a smaller order of magnitude than different spanning trees. Moreover, this assumption does not change significantly the ordering of spanning trees based on the associated line losses. 
\vspace{0.1cm}
\begin{assumption}
The impact of line losses on line flows is negligible relative to the power demands at the buses of the network.
\end{assumption}
This assumption implies that the power flow on each line $e\in \Eset$ is almost equal to the total demand of the buses that are receiving power through that line. Specifically, for a given spanning tree, if we denote the set of successors of an edge $e\in \Eset$ by $c_e$, then we have:
\begin{equation}
\label{eq:successor}
P_e=\sum_{i\in c_e} p_i \quad \text{and} \quad Q_e=\sum_{i\in c_e} q_i,
\end{equation}
where $p_i$ and $q_i$ are the active and reactive power demands at bus $i$. Note that by $P_e$ we mean the power flowing on line $e$ in the direction from the root of the tree to the leaves (parent to child). In Section \ref{sec:simulation}, we verify the validity of these assumptions in detail.

If we denote the loss of line $e=\{i,j\}\in \Eset$ by $L_{e}$, then we have:
 \begin{equation}
 \label{eq:loss}
 L_{e}=R_e \times \abs{I_{ij}}^2.
 \end{equation}
 In addition, in the relaxed model we have:
  \begin{equation}
  \label{eq:power}
  \abs{V_{i}}^2 \abs{I_{ij}}^2=P_{ij}^2+Q_{ij}^2.
 \end{equation}
 Combining (\ref{eq:successor}), (\ref{eq:loss}), and (\ref{eq:power}) with Assumption \ref{assumption:drop} implies that:
\begin{equation*}
L_e=R_e\left[ \left( \sum_{i\in c_e} p_i\right)^2+\left( \sum_{i\in c_e} q_i\right)^2\right].
\end{equation*}
 Given a spanning tree (ST) we can sum up the line losses $L_e$ over all the edges of the tree to find the total loss. Thus, the optimal reconfiguration problem with the goal of loss minimization can be written as the following optimization problem:
\begin{equation*}
\min_{ST} \quad \sum_{e\in ST} R_e\left[ \left( \sum_{i\in c_e} p_i\right)^2+\left( \sum_{i\in c_e} q_i\right)^2\right],
\tag{P1}
\label{eq:p1}
\end{equation*}
where the minimization is over all the spanning trees of $\mathcal{G}(\Nset,\Eset)$.

\section{Hardness Result} \label{sec:hardness}
In this section we prove that the DNR problem is strongly NP-hard in general by a reduction from the \tp\ problem \cite{gary1979computers}. A computational hardness result is more powerful when it is derived for a more restricted setting---since the hardness implication then holds for any generalization of the setting. Here we derive a hardness result for the special case of unit demands, where the objective function of the optimization problem (\ref{eq:p1}) reduces to a simpler function that is just counting the number of successors. In particular we make the following assumptions:
\begin{align*}
R_e=1&\quad \quad \forall e\in \Eset,\\
p_i=1 &\quad \quad \forall i\in \Nset\backslash\{0\},\\
q_i=0 &\quad \quad \forall i\in \Nset\backslash\{0\}.
\end{align*}
Under these assumptions, the optimization problem (\ref{eq:p1}) reduces to the following problem: 
\begin{equation*}
\min_{ST} \quad \sum_{e\in ST} (\text{number of successors of $e$ in $ST$})^2.
\tag{P2}
\label{eq:p2}
\end{equation*}
Although these assumptions may not be realistic, they transform the problem into an explicit combinatorial problem (without any power flow variable or parameter) and help us to analyze the computational complexity of the reconfiguration problem. The resulting complexity applies then to more general and realistic settings, as mentioned above.

We show that even the unit-demand case is strongly NP-hard. We prove this by a reduction from the \tp~problem defined as follows.
\vspace{0.1cm}
\begin{definition}
{\bf (\tp)} In the \tp\ problem we have a multiset of $k=3m$ integers summing to $mB$ with each integer strictly between $B/4$ and $B/2$. The task is to partition these numbers into $m$ triplets each with a sum of $B$.
\end{definition} 
It is well known \cite{gary1979computers} that in the \tp\ problem, deciding whether a given multiset can be partitioned into balanced triplets or not, is strongly NP-complete, i.e., it is NP-complete even if the numbers are bounded by a polynomial in the length of the input.

\begin{figure}[t] 
   \centering
   \includegraphics[width=2.65in]{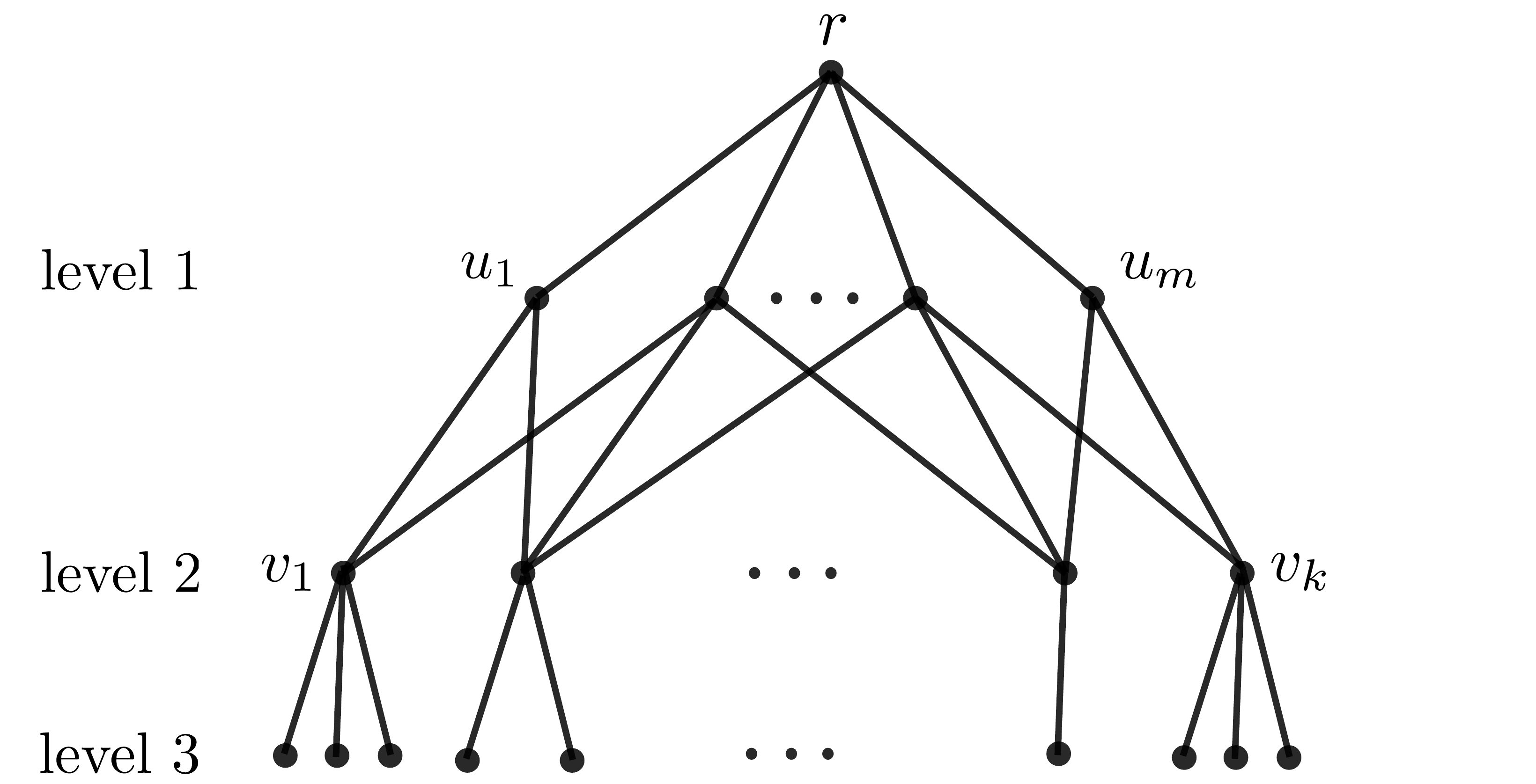} 
   \caption{Polynomial reduction.}
   \label{fig:reduction}
\end{figure}

\vspace{0.1cm}
\begin{theorem}[{\bf Hardness result}]
\label{theorem:hardness}
Distribution network reconfiguration problem (\ref{eq:p1}) is strongly NP-hard.
\begin{proof}
We propose a polynomial reduction from the \tp\ problem to the unit-demand case of reconfiguration problem (\ref{eq:p2}). Given an instance of the \tp~problem we build an instance of the reconfiguration problem such that the optimal spanning tree reveals the answer to the \tp\ problem (if it exists). Given $k=3m$ integers $\{a_1,a_2,...,a_{k}\}$, we construct a network as shown in Fig.~\ref{fig:reduction}. There is a root $r$, $m$ nodes $u_1,...,u_m$ connected to the root, $k=3m$ nodes $v_1,...,v_k$ each connected to all of $u_i$'s (thus $v_i$'s and $u_j$'s form a complete bipartite graph) and for each $v_i$ we have $a_i-1$ nodes connected to it. 

Lemma~\ref{lemma:rootedges} below proves that all the lines between the root $r$ and the $u_j$'s are part of the optimal tree. Moreover, all the lines between levels two and three appear in every spanning tree, so the only choices are on the lines between levels one and two. In particular, we have to connect each $v_i$ to exactly one $u_j$, i.e., make one of $m$ choices.

When we connect node $v_i$ to node $u_j$, the corresponding edge gets a cost of $a_i^2$, since there are $a_i-1$ nodes in level three and the edge has $a_i$ successors including $v_i$. This cost is independent of the choice of $u_j$, so the total cost for the edges between levels one and two is the same for all the spanning trees. The cost to be minimized is thus the total cost of the edges between root $r$ and the $u_j$'s. 

Let $S_j$ be the set of indices of the children of $u_j$, i.e.,
$$S_j=\{i\mid (u_j,v_i)\in \text{Tree}\},$$
then the cost related to edge $(r,u_j)$ is:
$$C_{(r,u_j)}=\left(1+\sum_{i \in S_j}a_i\right)^2,$$
where $1$ counts for the node $u_j$ itself. Now the total cost of the spanning tree is:
\begin{multline}
\label{eq:cost}
C=\sum_{j=1}^m C_{(r,u_j)}+\sum_{i=1}^{3m}a_i^2+\sum_{i=1}^{3m}(a_i-1)\\=\sum_{j=1}^m \left(1+\sum_{i \in S_j}a_i\right)^2+\sum_{i=1}^{3m}a_i^2+(mB-3m).
\end{multline}
As mentioned earlier, the second and third terms are constants since they are independent of the choice of the spanning tree.
Using the fact that the $S_j$'s are disjoint and $\cup_{j=1}^m S_j=\{1,...,k\}$, we also have:
\begin{multline}
\label{eq:constant}
\sum_{j=1}^m \left(1+\sum_{i \in S_j}a_i\right)=m+\sum_{j=1}^m \sum_{i \in S_j}a_i\\=m+\sum_{i=1}^{3m}a_i=m+mB=m(B+1).
\end{multline}

\begin{lemma}
\label{lemma:quadratic}
The minimum of $\sum_{i=1}^n x_i^2$ given that $\sum_{i=1}^n x_i=C$ for a constant $C\in \mathbb{R}$ is achieved when $x_i=C/n$ for all $1\leq i\leq n$.
\end{lemma}
By Lemma \ref{lemma:quadratic} and (\ref{eq:constant}), the minimum possible cost of (\ref{eq:cost}) is obtained when:
$$1+\sum_{i \in S_j}a_i=B+1 \quad \quad \forall j\in \{1,...,m\},$$
and the optimal value is:
$$C_{min}=m(B+1)^2+\sum_{i=1}^{3m}a_i^2+(mB-3m).$$
Note that this optimal cost is achieved when $a_i$'s are partitioned into $m$ subsets with sum $B$, but there is no restriction on the size of $S_j$'s. This means that node $u_j$ can have any number of $v_i$'s connected to it, while in the \tp\ problem we want to partition the $a_i$'s into $m$ triplets. The property $B/4<a_i<B/2$ ensures that this minimum can only be achieved when $|S_j|=3$ for all $j$. If for any $j'$ we have $|S_{j'}|>3$, then we get:
\small
$$\sum_{i \in S_{j'}}a_i>4 \times \frac{B}{4}=B,$$
\normalsize
and the partition cannot be balanced. Similarly if $|S_{j'}|<3$, then we get:
\small
$$\sum_{i \in S_{j'}}a_i<2 \times \frac{B}{2}=B.$$
\normalsize
In conclusion, the algorithm for the unit-demand case finds the tree corresponding to the \tp\ answer (if it exists), and if it outputs some unbalanced tree, this means that the \tp\ does not exist.
If each $a_i$ is bounded by a polynomial in $k$, the constructed network has polynomial number of nodes, hence any polynomial time algorithm for the unit-demand case provides a pseudo-polynomial time algorithm for the \tp\ problem which is not possible unless $P=NP$.
\end{proof}
\end{theorem}
Lemma~\ref{lemma:rootedges} proves the only remaining part of the hardness proof.
\vspace{0.1cm}
\begin{lemma}
\label{lemma:rootedges}
With uniform line resistances ($R_e=R,\forall e\in\Eset$), the optimal tree includes all the edges adjacent to the root.
\begin{proof}
We prove this by contradiction. Assume that we have an optimal tree which does not choose edge $(r,u)$ as shown in Fig.~\ref{fig:proof_lemma} on the left. 
\begin{figure}[t] 
   \centering
   \includegraphics[width=2.4in]{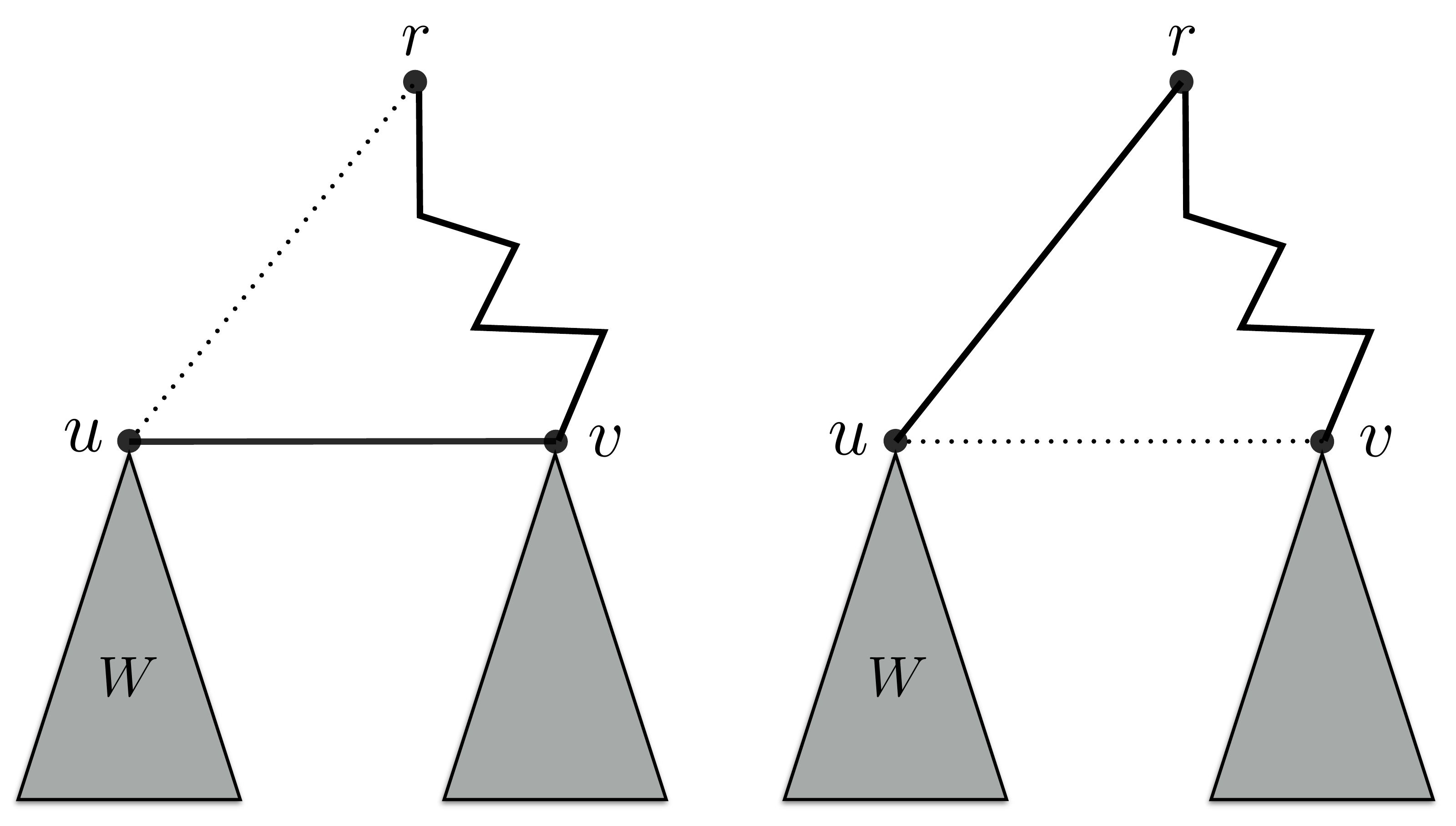} 
   \caption{Proof of Lemma \ref{lemma:rootedges}.}
   \label{fig:proof_lemma}
\end{figure}
Let $W$ be the total load weight of subtree connected to $u$ (including $u$). Since we have a tree, this subtree is connected to the root through another node $v$. Node $v$ may have other children and also may be connected to the root via one or more edges. Now we claim that this tree cannot be optimal since we can exchange edge $(u,v)$ with edge $(r,u)$ and improve the objective value as shown in the right tree. To see this, note that both edges $(u,v)$ in the left tree and $(r,u)$ in the right tree have costs $RW^2$, but the exchange of $(u,v)$ with $(r,u)$ decreases the load on all the edges of the path from $r$ to $v$ by $W$, hence decreasing the total cost. This contradicts the optimality of the first tree.
\end{proof}
\end{lemma}
Note that the unit-demand case is a special case of Lemma~\ref{lemma:rootedges}.
\section{Supermodular Structure}\label{sec:supermodular}
In the previous section we showed that DNR is strongly NP-hard, but if we find some additional structure such as submodularity or supermodularity in the problem, we may be able to provide approximation algorithms, which provides a rigorous worst-case performance guarantee. Here we first define this structure in addition to some required background about matroid constraints, and then we show that the DNR problem has this structure.
\subsection{Submodularity and Supermodularity} \label{subsec:submodularity}
Let $V$ be a finite set, called the ground set. We use $2^V$ to denote the set of all subsets of $V$, called the power set. A set function $f:2^V\mapsto \mathbb{R}$ is submodular if it has the \emph{diminishing returns} property, namely adding an element to a bigger set is less valuable than adding it to a smaller set.
\begin{definition}[{\bf Submodularity}]
A set function $f:2^V\mapsto \mathbb{R}$ with a ground set $V$ is submodular if:
$$f(X\cup \{u\})-f(X)\geq f(Y\cup \{u\})-f(Y),$$
for every $X\subseteq Y \subseteq V$, $u\in V\backslash Y$.
\end{definition}
Function $f$ is said to be supermodular if $-f$ is submodular (or the above inequality holds in the other direction). Supermodularity captures an increasing returns property. A function is said to be modular if it is both submodular and supermodular.
\vspace{0.1cm}
\begin{definition}[{\bf Monotonicity}]
A set function $f:2^V\mapsto \mathbb{R}$ is said to be monotone increasing if $f(X)\leq f(Y)$ for any $X\subseteq Y\subseteq V$.
\end{definition}
\vspace{0.1cm}
\begin{definition}[{\bf Matroid} \cite{schrijver2002combinatorial}]
Let $V$ be a finite set, and let $\mathcal{I}$ be a collection of subsets of $V$. The pair $\mathcal{M}=(V,\mathcal{I})$ is a matroid if the following conditions hold: (1) If $B\in \mathcal{I}$, then $A\in \mathcal{I}$ for all $A\subseteq B$, (2) If $A,B\in \mathcal{I}$ and $\abs{A}<\abs{B}$, then there exists $v\in B\backslash A$ such that $A\cup \{v\}\in \mathcal{I}$.
\end{definition}
A set $A\in \mathcal{I}$ is called an independent set. The collection $\mathcal{I}$ is called the set of independent sets of the matroid $\mathcal{M}$. A maximal independent set (an independent set that has maximum size) is a base of the matroid. It is easy to show that all the bases of a matroid have the same number of elements. 


\subsection{Set Function Formulation of DNR}
Considering the formulation of DNR (\ref{eq:p1}), the optimization problem is over all the spanning trees of the original graph. We would like to encode the two properties of ``being a tree'' and ``touching all the vertices of the graph'' into a set of constraints. In order to do this, we need to define a set of variables as follows. These variables also help to determine the successors of an edge in any arbitrary tree.
\begin{itemize}
\item For any edge $e\in \Eset$ we define a variable $x_{e}$ that indicates if that edge is included in the tree or not (the number of variables is equal to the number of lines in the distribution network).
\item Corresponding to any variable $x_{e}$, where $e=\{i,j\}$, we also define $y_{ij}^k$ and $y_{ji}^k$ for all $k\in \Nset$, which indicate the position of node $k$ compared to edge $e=\{i,j\}$. If there is a simple path from $i$ to $k$ including $\{i,j\}$, then $y_{ij}^k=1$ and if there is a simple path from $j$ to $k$ including $\{i,j\}$, then $y_{ji}^k=1$. In other words, $y_{ij}^k=1$ means that edge $\{i,j\}$ is chosen and $j$ is on the path from $i$ to $k$. If $x_{e}=0$, then both $y_{ij}^k$ and $y_{ji}^k$ are zero.
\end{itemize}
The following theorem, inspired by the integer programming formulation for the minimum spanning tree problem \cite{martin1986sharp,raghavan1994formulations}, explains how we use these variables to characterize the spanning trees explicitly.
\vspace{0.1cm}
\begin{theorem}[{\bf Feasible set characterization}]
There is a one-to-one correspondence between the spanning trees of $\mathcal{G}(\Nset,\Eset)$ and the feasible set specified by the following set of constraints:
\begin{align}
\sum_{e\in \Eset}x_{e}=n-1&\label{eq:c1} \\
y_{ij}^k+y_{ji}^k=x_{e} &\quad\quad \forall e=\{i,j\}\in \Eset, \forall k\in \Nset \label{eq:c2}\\
x_{e}+\sum_{k\neq i,j}y_{ik}^j =1&\quad\quad \forall i,j\in \Nset:e=\{i,j\}\in \Eset \label{eq:c3}\\
x_{e},y_{ij}^k,y_{ji}^k\in \{0,1\}&\quad\quad \forall e=\{i,j\}\in \Eset, \forall k\in \Nset \label{eq:c4}
\end{align}
\end{theorem}

If we write the total loss as a function of the binary variables above, we end up with an integer program formulation of (\ref{eq:p1}). 
For a given spanning tree $T$ (equivalently, a feasible set of values for the binary variables), and an edge $e=\{i,j\}\in T$, the variables $y_{ij}^k$ and $y_{ji}^k$ induce a partition of the vertices $\Nset$ into two sets which are exactly the two connected components of the tree obtained by removing $\{i,j\}$. The set that does not include the root (assuming vertex 0 is the root), is the set of successors of $e$ in $T$. In other words, if $y_{ji}^0=1$, and $c_e$ is the set of its successors, then we have:
\begin{equation*}
c_e=
\{k\in \Nset:y_{ij}^k=1\}.
\end{equation*}
Note that $y_{ji}^0=1$ is not an additional assumption, since the edges are not directed, and hence for the edges in the tree, one of the pairs $(i,j)$ or $(j,i)$ satisfies this condition.\\
Using this new description of successors, we can rewrite the objective function as:
\small
\begin{multline}
\label{eq:obj}
\sum_{e\in ST} R_e\left[ \left( \sum_{i\in c_e} p_i\right)^2+\left( \sum_{i\in c_e} q_i\right)^2\right]=\\
\sum_{i,j:\{i,j\}\in \Eset} R_{ij}y_{ji}^0 \left[ \left( \sum_{k\in \Nset} y_{ij}^k p_k\right)^2+ \left( \sum_{k\in \Nset} y_{ij}^k q_k\right)^2\right],
\end{multline}
\normalsize
where the inner summations are over all nodes, but the $y_{ij}^k$'s guarantee that we only count the successors, and the term $y_{ji}^0$ outside guarantees that we calculate each edge of the tree exactly once and in the correct direction with respect to the root.

So, the following optimization problem is equivalent to (\ref{eq:p1}):\footnote{we use $\sum_{\{i,j\}\in \Eset}$ instead of $\sum_{i,j\in \Nset:\{i,j\}\in \Eset}$ for simplicity.}
\begin{equation}
\small
\begin{split}
\min & \quad \sum_{\{i,j\}\in \Eset} R_{ij}y_{ji}^0 \left[ \left( \sum_{k\in \Nset} y_{ij}^k p_k\right)^2+ \left( \sum_{k\in \Nset} y_{ij}^k q_k\right)^2\right] \\
s.t. & \quad (\ref{eq:c1}),(\ref{eq:c2}),(\ref{eq:c3}),(\ref{eq:c4}).
\end{split}
\label{eq:p3}
\tag{P3}
\end{equation}

Now we show that (\ref{eq:p3}) is equivalent to a supermodular minimization problem with a single matroid basis constraint. The objective function (\ref{eq:obj}) is not supermodular over $\Eset$, but we create a similar set function that is supermodular and is equal to (\ref{eq:obj}) when constraints (\ref{eq:c1}--\ref{eq:c4}) hold (i.e., for spanning trees). A corollary to the following theorem shows that the feasible set in (\ref{eq:p3}) is indeed a matroid basis constraint. 
\vspace{0.1cm}
\begin{theorem}[{\bf Cycle Matroid} \cite{schrijver2002combinatorial}]
\label{theorem:cycle_matroid}
Let $\mathcal{G}(\Nset,\Eset)$ be an undirected graph. Define the set $T$ to be the collection of all subsets of $\Eset$ that form a forest (i.e., the subset is acyclic). In other words, $A\in T$ iff $A\subseteq \Eset$ and edges in $A$ do not form a cycle. Then $\mathcal{M}=(\Eset,T)$ is a matroid called the cycle matroid of graph $\mathcal{G}$ (also known as graphic matroid).
\end{theorem}
\vspace{0.1cm}
\begin{corollary}
\label{cor:basis}
Assuming that graph $\mathcal{G}$ is connected, the bases of the cycle matroid $\mathcal{M}$ are the spanning trees of $\mathcal{G}$, which all have cardinality $\abs{\Nset}-1$. Therefore, constraints (\ref{eq:c1}--\ref{eq:c4}) are equivalent to a single matroid basis constraint on $\Eset$.
\end{corollary}
Now we introduce the supermodular set function over $\Eset$. For any $A\subseteq \Eset$, we define:
\begin{equation}
\small
\label{eq:sup_obj}
f(A)=\sum_{\{i,j\}\in \Eset} R_{ij}z_{ji}^0 \left[ \left( \sum_{k\in \Nset} z_{ij}^k p_k\right)^2+ \left( \sum_{k\in \Nset} z_{ij}^k q_k\right)^2\right]
\end{equation}
The only difference between (\ref{eq:obj}) and (\ref{eq:sup_obj}) is that we replaced the $y_{ij}^k$'s with $z_{ij}^k$'s, and $z_{ij}^k$ is defined similar to $y_{ij}^k$ except that it can be any non-negative integer (compared to $0,1$) and it counts the number of paths in $A$ starting with $\{i,j\}$ and going to $k$. Clearly, for spanning trees there cannot be more than one path between any arbitrary pair of vertices, therefore $z_{ij}^k=y_{ij}^k$ and this implies the equality of (\ref{eq:obj}) and (\ref{eq:sup_obj}) when constraints (\ref{eq:c1}--\ref{eq:c4}) hold.
\vspace{0.1cm}
\begin{theorem}[{\bf Supermodularity}]
\label{theorem:supermod}
Objective function (\ref{eq:sup_obj}) is a supermodular set function over $\Eset$, provided that the $p_i$'s and $q_i$'s are non-negative.
\begin{proof}
The sum of supermodular set functions is supermodular, so we only need to prove the supermodularity for a fixed edge $\{i,j\}\in \Eset$. We can also drop positive constants like $R_{ij}$. Define $f_{ij}(A)$ and $f_{ij}'(A)$ as follows:
\small
\begin{equation*}
\label{eq:fij}
f_{ij}(A)=z_{ji}^0 \left( \sum_{k\in \Nset} z_{ij}^k p_k\right)^2, f_{ij}'(A)=z_{ji}^0 \left( \sum_{k\in \Nset} z_{ij}^k q_k\right)^2
\end{equation*}
\normalsize
We now prove that $f_{ij}(A)$ is supermodular. A similar proof works for $f_{ij}'(A)$. We want to show that:
\begin{equation}
\label{eq:sup_want}
f_{ij}(A\cup \{e\})-f_{ij}(A)\leq f_{ij}(B\cup \{e\})-f_{ij}(B),
\end{equation}
for every $A\subseteq B \subseteq \Eset$, $e\in \Eset$, and $e\not\in B$. For any $k$, let $a_{ij}^k$ be the change in $z_{ij}^k$ when we add $e$ to $A$, i.e.:
$$a_{ij}^k=z_{ij}^k(A\cup \{e\})-z_{ij}^k(A),$$
where $z_{ij}^k(A)$ is just $z_{ij}^k$, calculated based on the edges in $A$. Similarly, let $b_{ij}^k$ be defined for $B$ and $B\cup \{e\}$. We have $a_{ij}^k\leq b_{ij}^k$, because any new path in $A$ created by adding $e$ is also a new path in $B$. Another fact is that $z_{ij}^k(A)\leq z_{ij}^k(B)$, because adding more edges cannot decrease the number of paths between any pair of vertices (i.e., $f(A)$ is a monotone increasing function). Now we prove (\ref{eq:sup_want}):
\small
\begin{align}
&f_{ij}(B\cup \{e\})-f_{ij}(B)\\
\begin{split}
=&z_{ji}^0(B\cup \{e\}) \left( \sum_{k\in \Nset} z_{ij}^k(B\cup \{e\}) p_k\right)^2 \label{eq:step1}\\
&-z_{ji}^0(B) \left( \sum_{k\in \Nset} z_{ij}^k(B) p_k\right)^2
\end{split}\\
\begin{split}
=&(z_{ji}^0(B)+b_{ji}^0) \left(\sum_{k\in \Nset} b_{ij}^k p_k+\sum_{k\in \Nset} z_{ij}^k(B) p_k\right)^2 \label{eq:step2}\\
&-z_{ji}^0(B) \left( \sum_{k\in \Nset} z_{ij}^k(B) p_k\right)^2
\end{split}\\
\begin{split}
\geq&(z_{ji}^0(A)+b_{ji}^0) \left(\sum_{k\in \Nset} b_{ij}^k p_k+\sum_{k\in \Nset} z_{ij}^k(B) p_k\right)^2 \label{eq:step2.5}\\
&-z_{ji}^0(A) \left( \sum_{k\in \Nset} z_{ij}^k(B) p_k\right)^2
\end{split}\\
\begin{split}
\geq &(z_{ji}^0(A)+b_{ji}^0) \left(\sum_{k\in \Nset} b_{ij}^k p_k+\sum_{k\in \Nset} z_{ij}^k(A) p_k\right)^2 \label{eq:step3}\\
&-z_{ji}^0(A) \left( \sum_{k\in \Nset} z_{ij}^k(A) p_k\right)^2
\end{split}\\
\begin{split}
\geq &(z_{ji}^0(A)+a_{ji}^0) \left(\sum_{k\in \Nset} a_{ij}^k p_k+\sum_{k\in \Nset} z_{ij}^k(A) p_k\right)^2 \label{eq:step4}\\
&-z_{ji}^0(A) \left( \sum_{k\in \Nset} z_{ij}^k(A) p_k\right)^2
\end{split}\\
=&f_{ij}(A\cup \{e\})-f_{ij}(A).
\end{align}
\normalsize
In (\ref{eq:step1}), (\ref{eq:step2}) we just applied the definitions of $f_{ij}$ and $b_{ij}^k$, respectively. In (\ref{eq:step2}), the aggregate coefficient of $z_{ji}^0(B)$ is positive, so using the fact that $z_{ji}^0(B)\geq z_{ji}^0(A)$, we get (\ref{eq:step2.5}). To get (\ref{eq:step3}), note that quadratic function $(x+\alpha)^2-x^2<(y+\alpha)^2-y^2$ for $x<y$ and fixed $\alpha>0$. Setting $\alpha=\sum_{k\in \Nset} b_{ij}^k p_k$, and the fact that $z_{ij}^k(A)\leq z_{ij}^k(B)$ implies (\ref{eq:step3}). Finally, (\ref{eq:step4}) is implied by $a_{ij}^k\leq b_{ij}^k$, which proves the supermodularity of $f_{ij}(A)$. We have
$$f(A)=\sum_{\{i,j\}\in \Eset} R_{ij}\bigg(f_{ij}(A)+f_{ij}'(A)\bigg),$$
therefore $f(A)$ is supermodular.
\end{proof}
\end{theorem}
\section{Algorithm and Performance Guarantee}\label{sec:alg}
\begin{table*}[ht]
    \centering
    \caption{Comparison of different DNR algorithms on the 33-bus network of Fig.~\ref{fig:33bus}}
    \label{table:exp}
    \small
    \begin{tabular}{| c | c | c | c |}
        \hline
         \rowcolor{lightgray} Algorithm						& Method 					& Open Lines 			& Loss (kW) \\
        \hline
             Proposed 									& Submodular Local Search	& $7,9,14,32,37$ 		& $139.552$\\
        \hline
            Morton and Mareels \cite{morton2000efficient} 		& Brute-Force 				& $7,9,14,32,37$ 		& $139.552$\\
        \hline
            Gomes \etal \cite{gomes2005new} 				& Greedy on Mesh Network 	& $7,9,14,32,37$ 		& $139.552$\\
        \hline
            Khodr and Martinez-Crespo \cite{khodr2009integral} 	& Benders Decomposition		& $7,9,14,32,37$ 		& $139.552$\\
        \hline
           Wu \etal \cite{wu2010study}						& Ant Colony				& $7,9,14,28,32$		& $139.976$\\
        \hline
            Shirmohammadi and Hong \cite{shirmo1989reconfig} 	& Optimal Current Pattern 	& $7,10,14,32,37$ 		& $140.279$\\
        \hline
            Baran and Wu \cite{baran1989network}\textsuperscript{\ref{notes:note1}} 			& Branch Exchange 			& $11,28,31,33,34$ 		& $146.832$\\
        \hline
        \multicolumn{2}{|c|}{Initial Configuration} 										& $33,34,35,36,37$		& $202.670$\\
        \hline
    \end{tabular}
\end{table*}
In the previous section we showed that the DNR problem (\ref{eq:p3}) is equivalent to a supermodular minimization problem subject to a single matroid basis constraint. Unless $P=NP$, it is not possible to approximate the minimum of a supermodular function within any factor \cite{mittal2013fptas}, in contrast with the related problem of maximizing a submodular function which admits a constant factor approximation algorithm \cite{lee2009non}. We adapt the approximation algorithm for the submodular maximization problem under matroid constraints, proposed by Lee~\etal \cite{lee2009non}, to solve the DNR problem, but we have to convert the supermodular function to a non-negative submodular function (by negating and shifting). This conversion affects the multiplicative approximation guarantee, as shown in Theorem \ref{theorem:approx}. The algorithm, which is based on local search, is described in Algorithm \ref{alg:localsearch}.
\begin{algorithm}[h]
\caption{Distribution Network Reconfiguration for Loss Minimization}
\label{alg:localsearch}
\begin{algorithmic}[1]
	\STATE \textbf{Input:} Configuration $\mathcal{G}(\Nset,\Eset)$, bus demands ($p_k,q_k$), line resistances ($R_{ij}$), $\epsilon$.
	\STATE \textbf{Output:} Spanning tree for minimizing the total loss.
	\STATE \textbf{Initialize} $T$ with an arbitrary spanning tree.
	\WHILE{1}
		\IF{there exist $e\in \Eset\backslash T$ and $e' \in T$ such that $(T\backslash\{e'\})\cup\{e\}$ is a spanning tree and $f((T\backslash\{e'\})\cup\{e\})<(1-\epsilon)f(T)$}
			\STATE $T\leftarrow (T\backslash\{e'\})\cup\{e\}$
		\ELSE
			\STATE break
		\ENDIF
	\ENDWHILE
	\RETURN $T$
\end{algorithmic}
\end{algorithm}

The algorithm starts with an arbitrary spanning tree $T$. Then at each iteration, it looks for two edges $e\in \Eset \backslash T$ and $e'\in T$ such that swapping those two edges makes another spanning tree with loss at most $(1-\epsilon)f(T)$. If such a pair exists, it updates $T$ and repeats the exchange process, otherwise the algorithm terminates and outputs the locally optimal spanning tree.
\vspace{0.1cm}
\begin{theorem}[{\bf Performance guarantee}]
\label{theorem:approx}
Let $T_{alg}$ be the output of Algorithm \ref{alg:localsearch}, and $T^*$ be the optimal spanning tree, i.e., $T^*=argmin\{f(T):T\subseteq \Eset, T \text{ is a spanning tree}\}$. Let $M=f(\Eset)$, which is an upper bound on $f(A)$ for all $A\subseteq \Eset$, then:
\begin{equation}
M-f(T_{alg})\geq \left(\frac{1}{6}-\epsilon\right)\big(M-f(T^*)\big).
\end{equation}
\begin{proof}
This is a corollary of \cite[Theorem~22]{lee2009non}, which provides a $(\frac{1}{6}-\epsilon)$-approximation algorithm for maximizing any non-negative submodular function over bases of a matroid $\mathcal{M}$.\footnote{That theorem requires $\mathcal{M}$ to have at least two disjoint bases. We can solve this (if necessary) by adding dummy edges with very high resistances (to make sure that the algorithm never selects them). Moreover, their algorithm performs another local search which allows deletion of elements, but that run yields the empty set in our case (due to the monotonicity), hence does not apply to the DNR problem.} Here we use $M-f(A)$ as the non-negative submodular function, and the spanning trees are the bases of the cycle matroid discussed in Theorem \ref{theorem:cycle_matroid}. 
\end{proof}
\end{theorem}
Even though Algorithm \ref{alg:localsearch} is based on the local search approximation algorithm for maximizing non-monotone submodular functions \cite{lee2009non}, it is equivalent to the branch exchange heuristic algorithm which has been used since the late 1980s \cite{baran1989network}. This establishes that Theorem \ref{theorem:approx} provides the first proof of a performance bound, and hence a performance guarantee for the branch exchange algorithm.
\section{Experiments}\label{sec:simulation}
Table \ref{table:exp} shows the results of our experiments on the 33-bus network of Fig.~\ref{fig:33bus}. The parameters of the network can be found in \cite{pso}. All the active and reactive power demands are positive for this network as assumed in Theorem \ref{theorem:supermod}. The simulations have been done by using the MATPOWER package in \textsc{Matlab} \cite{zimmerman2011matpower}. The results show that in this case, our submodular approach finds the globally optimal configuration, which was found in \cite{morton2000efficient} (by enumerating all 50751 spanning trees). In \cite{baran1989network}, 2 different approximate power flow methods with different accuracies have been used and we also believe that there are inconsistencies regarding the parameters of the network in the literature\footnote{\label{notes:note1}The resistance of the branch between bus 6 and bus 7 is $0.7114 \Omega$ in \cite{baran1989network}, but $1.7114 \Omega$ in \cite{pso}. We used the latter value in all our simulations.}; that is why results reported in \cite{baran1989network} differ from what we obtained by Algorithm~\ref{alg:localsearch}. Clearly, the output of the local search algorithms depends on the initialization. We used the initial configuration (Fig.~\ref{fig:33bus}) as the initial spanning tree in our simulation. Further, to check the robustness with respect to the initial tree, we repeated the simulations with $1000$ random initial trees, all of which ended with the same optimal solution.

\begin{figure}[t] 
   \centering
   \includegraphics[width=1\linewidth]{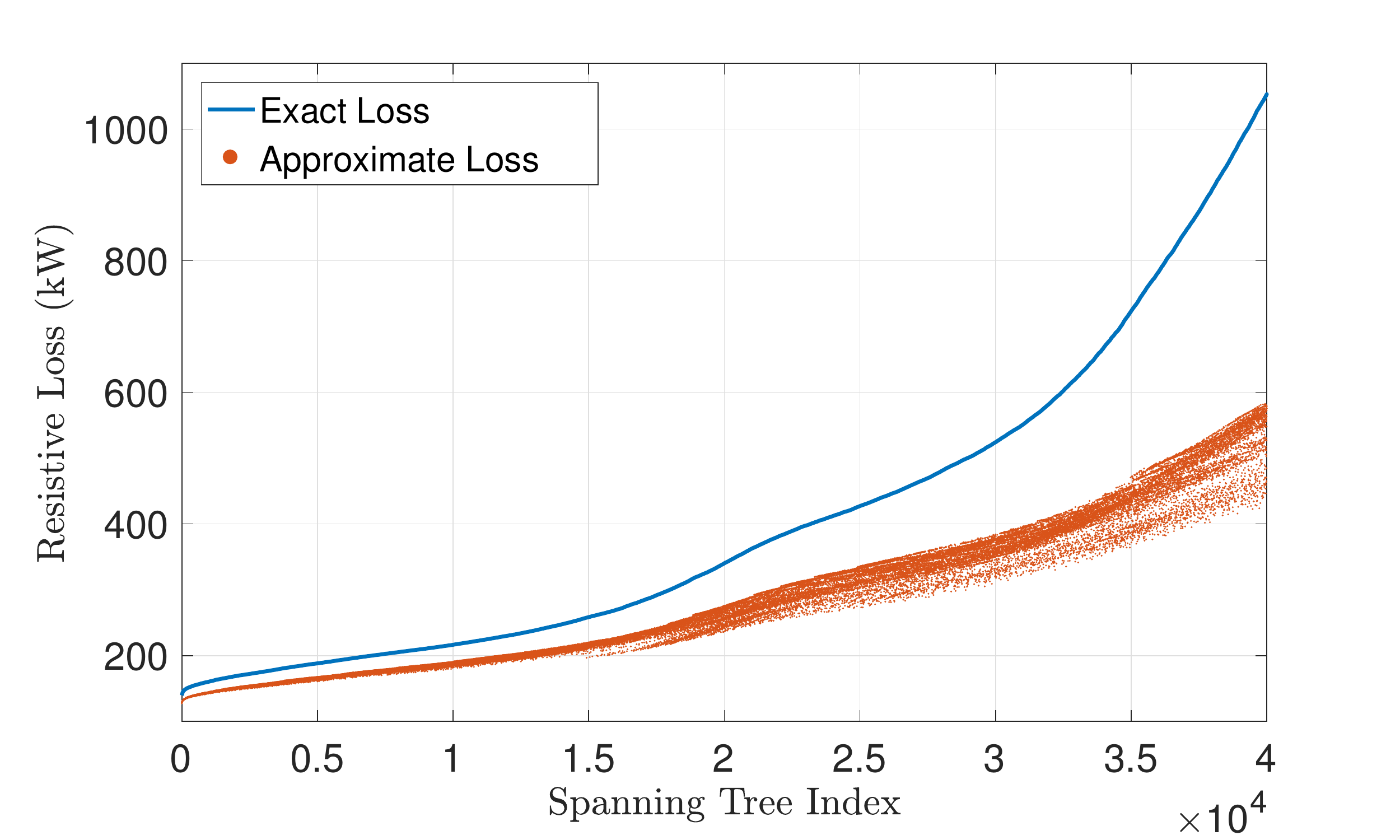} 
   \caption{Comparing losses from the simplified model with the exact values.}
   \label{fig:loss_comparison}
\end{figure}

In order to check the validity of our assumptions (see the problem formulation in Section \ref{sec:model}), we compare the losses of spanning trees as measured in (\ref{eq:p1}) with the exact losses obtained from MATPOWER. The result is shown in Fig.~\ref{fig:loss_comparison}. The blue line is the exact loss curve where the spanning trees are sorted in the order of increasing total loss. The red dots also show the loss for each tree obtained from the simplified model. We observe that the approximate loss is generally increasing, which means that it can be used in the local search algorithm. In fact performing the local search with either exact loss or approximate loss results in the same globally optimal tree, reported in Table \ref{table:exp}.
As expected, approximate loss estimates are less accurate for trees with higher losses, since the resistive losses approach the order of magnitude of load demands in such networks (hence contradicting Assumption $2$).
On the other hand, for trees with smaller losses (which are indeed the target of our optimization problem) the simplified loss approximates the exact loss very well.

\begin{figure}[t] 
   \centering
   \includegraphics[width=0.8\linewidth]{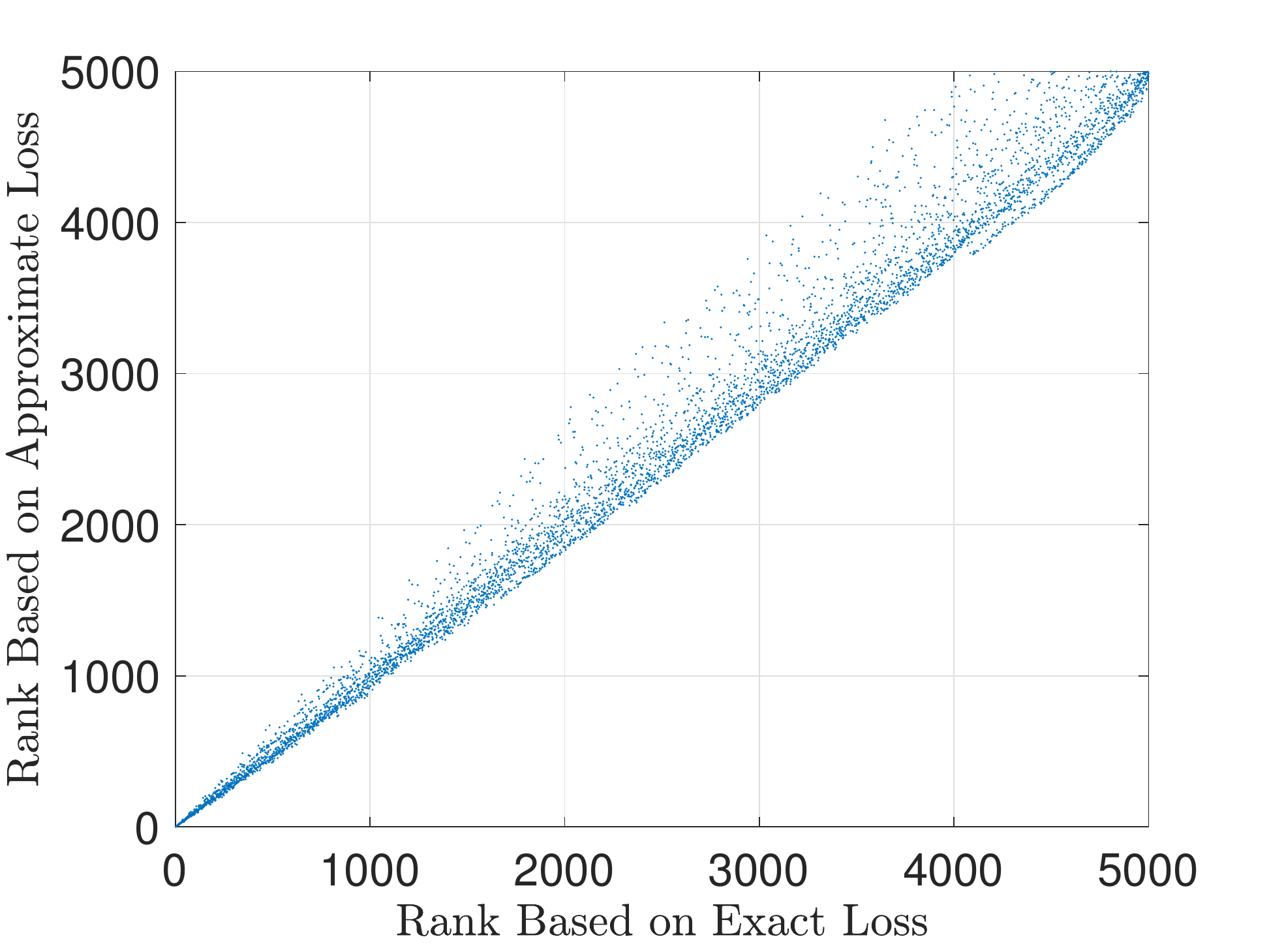} 
   \caption{Rankings based on exact and approximate losses for the best $5000$ spanning trees.}
   \label{fig:rank}
\end{figure}

Fig.~\ref{fig:rank} also compares the rank of the top 5000 spanning trees based on the exact and approximate losses. Ideally, we would like the simplified losses to preserve the rankings (which would result in a $y=x$ line in this plot). We observe that no single spanning tree faces a significant change in its ranking.
\section{Conclusion} \label{sec:conclusion}
\vspace{-0.1in}
In this paper, we studied the distribution network reconfiguration problem (DNR) for loss minimization through a submodular optimization approach. We proved that this problem is NP-hard even if the demands and the line resistances are all equal to one. We formulated this problem as a supermodular minimization problem subject to a matroid basis constraint. We then used the algorithm for maximizing non-monotone submodular functions under matroid constraints, to give a polynomial time algorithm for the DNR problem with a performance guarantee. The algorithm was equivalent to the branch exchange algorithm that was known previously, but for which no theoretical guarantees were available. By discovering a submodular structure in the problem, we pioneered the derivation of a performance bound on the branch exchange algorithm.

Although supermodular minimization cannot be approximated in general, there are approximation algorithms for the case when the supermodular function has bounded curvature (see \cite{il2001approximation,sviridenko2015optimal} for the definition of curvature and the approximation algorithms). The formulation studied in this paper does not have bounded curvature. One interesting question that arises is whether the DNR problem can be formulated as minimizing a supermodular function with bounded curvature. A positive determination would imply a multiplicative constant factor approximation (compared to Theorem \ref{theorem:approx} which includes the upper bound $M$) and would provide a significant improvement.


\balance
\section{References}
\vspace*{-6ex}
\footnotesize
\bibliographystyle{ieeetr} 
\bibliography{DSR}

\end{document}